\documentclass[11pt]{article}
\usepackage{amsmath,amssymb,bm,amsthm}
\usepackage{graphicx,xcolor}

\usepackage[style=numeric]{biblatex}
\renewbibmacro{in:}{}
\usepackage{array,url}
\usepackage{amsfonts,epstopdf}
\usepackage{booktabs} 
\usepackage{array}
\usepackage{caption} 
\usepackage[position=top]{subfig}
\usepackage[shortlabels]{enumitem}
\usepackage[bookmarks=false]{hyperref} 
\usepackage[stable]{footmisc}
\usepackage[capitalise,nameinlink,noabbrev]{cleveref}
\crefname{equation}{}{}
\let\ref\cref
\let\eqref\cref
\let\autoref\cref

\usepackage{indentfirst}
\usepackage{tikz}
\graphicspath{{./}{./figs/}{./figs/nets/}{./figs/efsketch/}}

\newcommand{\p}{\partial}

\newcommand{\rar}{\rightarrow}

\newcommand{\R}{\mathbb{R}}
\newcommand{\K}{\kappa}

\usepackage{algorithm}
\usepackage[noend]{algpseudocode}
\makeatletter
\def\BState{\State\hskip-\ALG@thistlm}
\makeatother
\usepackage{algorithmicx}

\title{\LARGE \bf
    Linking Machine Learning with Multiscale Numerics: 
    Data-Driven Discovery of Homogenized Equations}

\author{Hassan Arbabi, Judith E. Bunder, Giovanni Samaey, \\Anthony J. Roberts and Ioannis G. Kevrekidis
	\thanks{H.A. is with the Department of Mechanical Engineering,
		Massachusetts Institute of Technology, and Department of Chemical and Biomolecular Engineering at Johns Hopkins University.  J.E.B. and A.J.R. are with School of Mathematical Sciences, University of Adelaide. G.S. is with the Computer Science Department at KU Leuven.
		I.G.K. is with the Departments of Chemical and Biomolecular Engineering \& Applied Mathematics and Statistics at Johns Hopkins University.
		Correspondence: {\tt\small arbabiha@gmail.com, yannisk@jhu.edu}}
}

\hypersetup{
  colorlinks   = true, 
  urlcolor     = blue, 
  linkcolor    = blue, 
  citecolor   = magenta 
}

\newif\ifcomments

\commentsfalse

\newif\iftodo

\todotrue

\begin{document}

\maketitle

\begin{abstract}
The data-driven discovery of partial differential equations (PDEs) consistent with spatiotemporal data is experiencing a rebirth in machine learning research.  Training deep neural networks to learn such data-driven partial differential operators requires extensive spatiotemporal data. For learning coarse-scale PDEs from computational fine-scale simulation data, the training data collection process can be prohibitively expensive. We propose to transformatively facilitate this training data collection process by linking machine learning (here, neural networks) with modern multiscale scientific computation (here, equation-free numerics). These equation-free techniques operate over sparse collections of small, appropriately coupled, space-time subdomains (``patches"), parsimoniously producing the required macro-scale training data. Our illustrative example involves the discovery of effective homogenized equations in one and two dimensions, for problems with fine-scale material property variations.  The approach holds promise towards making the discovery of accurate, macro-scale effective materials PDE models possible by efficiently summarizing the physics embodied in ``the best" fine-scale simulation models available.

\end{abstract}


\clearpage

\section{Introduction }

Evolutionary models of materials dynamic behavior, in the form of Partial Differential Equations (PDEs), embodying conservation laws supplemented by appropriate closures, traditionally form the backbone of computational materials modeling. 
The requisite closures were initially mostly phenomenological, guided by experimentation, but in recent years such closures increasingly come from fine-scale, microscopic, possibly even \emph{ab initio} computations. 
One of the signal promises of machine learning in its early days, and increasingly more so today, is the discovery of such evolutionary PDEs from spatiotemporal data --- whether from physical observations or from computational observations obtained through fine-scale models.  

Yet materials problems are inherently multiscale: fine-scale models typically operate at the atomistic level; but the questions asked, and the answers required, are usually macroscopic in nature~\cite[e.g.][]{NASA2018}, and the models that describe  macroscopic materials physics are coarse-scale, effective PDEs in terms of macroscopic observables.

When the macroscopic Quantities-of-Interest (QoI), or the ``right variables", in terms of which an effective macroscopic PDE can be written, are known, there already exist data-driven, machine-learning-assisted techniques that will successfully approximate a PDE for the QoI evolution from macroscopic spatiotemporal observations. These techniques are based on the simple observation that \emph{the law of the PDE} is a functional relation (a constraint) between local space and time derivatives of the field(s) of interest: for simple diffusion, the  time derivative of the local density is a function of its second-order space derivative.
From detailed spatiotemporal data, one can thus acquire extensive sampling of instantiations of this relation: every space-time point in a movie provides a training data point for a machine learning model that can learn a PDE right-hand-side (i.e., a differential operator) consistent with the observed movie,  typically assuming the PDE is homogeneous in space-time. 
Such neural-network-based techniques were certainly proposed and implemented in the 1990s for ODEs \cite[e.g.][]{rico1992discrete,rico1995nonlinear} and PDEs \cite[e.g.][]{gonzalez1998identification}. Interestingly, several neural net architectural features (recurrent networks, convolutional networks) that proved broadly successful later, arose naturally in the work described in these references.
The field has exploded recently, and tools beyond deep neural networks (e.g., Gaussian Processes as well as Sparse Regression and Reservoir Computing) have been brought to the task   ~\cite{rudy2017data,raissi2018numerical,raissi2018hidden,pathak2018model,vlachas2018data,wang2019variational}. 
Gray box techniques (where part of the equation is known ~\cite[e.g.][]{psichogios1992hybrid,rico1994continuous,bar2019learning}) as well as machine-learning techniques for solving known PDEs  are increasingly successful \cite{lagaris1998artificial,raissi2019physics,kharazmi2019variational,meng2020ppinn}.

These approaches can learn approximate macroscopic PDEs based on data of macroscopic field evolution.
The data may come from macroscopic measurements of physical experiments,  or from coarse-grained observations of fine-scale simulations (e.g., running molecular dynamics, kinetic Monte Carlo or agent-based models, and observing the evolution of macroscopic, hydrodynamic level fields, such as densities, momenta or stresses).
We recently discussed the learning of coarse-grained PDEs based on fine-scale simulators: the ``inner", fine-scale simulator was a Lattice-Boltzmann (LB) model, while the ``outer", coarse-grained fields were two concentration fields for a system of coupled reaction-diffusion equations~
\cite{lee2020coarse}. 
The fine-scale LB simulation was performed ``for all space and all times" over the domain of interest, which could lead to considerable computational effort.
The goal of this paper it to alleviate, as much as possible, the extensive computational needs for the collection of the requisite macro-scale training data.
Here we work on coarse-graining a class of problems different than the more ``inner atomistic simulation" class illustrated by the LB example: we are interested in discovering coarse effective PDEs for materials that have micro-scale variations in their properties. In this case the ``inner simulator" is the PDE that resolves that complete fine-scale material properties, while the coarse-grained PDE is the so-called ``homogenized PDE" --- a PDE that governs the effective, long-wavelength material response. One may loosely think of this as a PDE for the locally averaged (over the fine-scale variations) material response. 
In applied mathematics, this problem is the purview of \emph{homogenization theory} \cite{bensoussan2011asymptotic} (see also \cite{Bunder2018a}): using asymptotic techniques, and under specific assumptions, one can derive a closed-form analytical expression for the homogenized PDE. The homogenized PDE solution is the part of the original solution that does not exhibit micro-scale variation. Therefore, the homogenized PDE is the right \emph{effective} coarse-scale PDE for materials with microscopic heterogeneity. However, the analytical derivation of homogenized PDEs is typically based
on assumptions and/or asymptotic arguments that may not be valid for many practical problems.
In such (frequent!) cases, the data-driven discovery of a useful approximate effective PDE, 
may be the only practical alternative. 
Our plan in this paper is to mitigate this ``training data collection" computation expense by taking advantage of modern multiscale scientific computing approaches (in particular, equation-free techniques) that are based on precisely the same assumption: that an effective equation exists, yet it is not available in closed form.

Equation-free multiscale computation comprises a suite of scientific computing techniques for solving macroscopic evolution equations based on parsimonious use of fine-scale simulators without ever deriving the macroscopic equations in closed form \cite{kevrekidis2003equation,samaey2009equation}.
Equation-free algorithms like coarse projective integration, gap-tooth and patch dynamics, as well as coarse bifurcation and stability analysis, are templated on established numerical techniques; the new element is that the quantities required for scientific computation (time derivatives, action of Jacobians) are not obtained from closed-form function evaluations, but rather from brief bursts of well-designed fine-scale simulation over short time intervals and small spatial domains.
It is precisely this ``closure on demand" property of equation-free computation that will help us collect the macroscopic data required for macroscopic PDE learning \slash approximation with significantly less computational effort. 
Such parsimonious design of micro-scale computations to collect macro-information adequate for learning, is the main ingredient of our work.

The rest of this paper is structured as follows. 
We start in \cref{sec_homeq} by describing homogenization problems arising from the modeling of materials with micro-scale property variation.
The equation-free methods we use are described in \cref{sec_eqfree}; they help us collect large-space\slash long-time data at the homogenized level from minimal (small-space\slash short-time) finely resolved direct material simulation. We then use the collected macro-data to learn the effective, macro-scale, homogenized equation through the use of two neural network architecture variations.
\cref{sec_learning} describes our approach to data-driven PDE learning and the neural network architectures.
In \cref{sec_results}, we apply our framework to two illustrative problems.  In each problem we compare the data-driven effective model we identify with analytically obtainable homogenized equations, thus establishing the utility but also highlighting some limitations of our approach. 
We conclude in \cref{sec_conclusion} with a brief discussion of possible extensions and applications of our approach and some thoughts on the important issue of data-driven coarse variable selection for problems where even the appropriate QoI must be discovered from the data. 
A python implementation of our framework and the illustrative examples in \cref{sec_results} is available at  \url{https://github.com/arbabiha/homogenization_via_ML}.

\clearpage
\section{Homogenization problem}\label{sec_homeq}
In many physical or engineering modeling problems, the properties of the material under investigation have variations (i.e., heterogeneity) at small spatial scales, while we are interested in predicting the response of the material at much larger scales. 
The \emph{homogenization technique} is an asymptotic applied mathematics method that, under certain technical assumptions, derives (approximate) effective equations for the material response at large length scales. Here, we briefly review this technique through the example of diffusion in a material with heterogeneous diffusivity---comprehensive discussions are given elsewhere~\cite{bensoussan2011asymptotic, pavliotis2008multiscale}.

Consider the unsteady diffusion problem for~\(u(x,t)\) in the one-dimensional spatial domain~$[0,1]$ with homogeneous Dirichlet boundary conditions:
\begin{align}\label{eq_detailed}
\p_t u(x,t)&=\p_x\left[ a\big(\frac{x}{\epsilon}\big) \p_x u(x,t)\right] , \notag\\
 u(x,0)&=u^0(x) \in L^2([0,1]),\\
 u(0,t)&=u(1,t)=0 , \notag
\end{align}
where the diffusivity~$a(\cdot)$ is space-dependent, and periodic in its argument. We are interested in the case of $\epsilon \ll1$, that is, the diffusivity has very rapid spatial oscillations. This is \emph{the detailed problem}.
Since~$\epsilon$ is very small compared to the length scales of interest (e.g., the size of the domain), the goal of homogenization is to find an effective \emph{homogenized PDE} that has no explicit dependence on~$\epsilon$, but its solution usefully approximates the solution to the detailed problem in its large-scale features. The classical approach to derive such a PDE is to expand the solution of the detailed problem in powers of~$\epsilon$,
\begin{align}\label{ep_expansion}
u(x,t) = u_h(x,t) + \sum_{j=1}^\infty \epsilon^j  u_j(x,x/\epsilon,t),
\end{align}
where the functions $u_j,~j=1,2,\ldots$, are periodic in their second argument. The first term in the expansion,~$u_h(x,t)$, is the homogenized solution and, to errors~\(O(\epsilon^2)\), it satisfies the homogenized PDE
\begin{align}\label{eq_hom}
\p_t u_h(x,t)&=\p_x\left[ a^* \p_x u_h(x,t)\right] , \notag\\
 u_h(x,0)&=u^0(x) \in L^2([0,1]),\\
 u_h(0,t)&=u_h(1,t)=0 , \notag
\end{align}
where \emph{the effective diffusivity constant} is 
\begin{align}\label{eq_effdiff}
a^*=\int_{0}^{1}a(y)\left[1-\frac{d}{dy}\chi(y)\right]dy.
\end{align}
In the above equation, $\chi(y)$ is the periodic solution to the so-called cell problem,
\begin{align}\label{eq_cell}
\frac{d}{dy}\left[a(y) \frac{d}{dy}\chi \right]=\frac{d}{dy}a(y),
\end{align}
which is usually accompanied by a normalization condition,
\begin{align}\label{eq_cellnormal}
\int_0^1\chi(y)dy=0.
\end{align}
It is known  that as $\epsilon\rar 0$, the homogenized solution~$u_h$ converges to the solution of the detailed problem~\cite{bensoussan2011asymptotic}, and that at finite~\(\epsilon\) one can systematically find corrections to~\eqref{eq_hom} for~\(u_h\)~\cite{Roberts2013a}.  Unlike the detailed field~$u$, the homogenized solution does not have an explicit dependence on the  fast-varying variable~$x/\epsilon$, and hence it represents the large-scale component of the detailed solution. When~$\epsilon$ is small enough, one might consider approximating the homogenized solution by applying a low-pass filter to the detailed solution (i.e., filtering out variations with wavelengths of order~$\epsilon$ or smaller). In this paper, we take the coarse scale (low-pass filtered, locally averaged)  density to be our effective macroscopic variable---which we expect approximates~$u_h$---and propose a framework to discover its evolution PDE from data.

Application of the rigorous analytical homogenization process, outlined above, to PDEs suffers from limitations because its validity requires somewhat stringent mathematical assumptions and asymptotic limits. For example, in the general case of non-periodic heterogeneity, no simple effective diffusivity formula can be derived; indeed, no systematic approach for deriving it yet exists~\cite{pavliotis2008multiscale}.  
On the other hand, direct numerical simulation of heterogeneous materials requires a computational power scaling as at least~$(1/\epsilon)^{d}$, with~$d$ being the spatial dimension, which makes such problems computationally prohibitive.

\section{Equation-free approach for homogenization} \label{sec_eqfree}
The equation-free approach, introduced by Kevrekidis et al.~\cite{kevrekidis2003equation}, when applied to homogenization problems~\cite{runborg2002effective, samaey2009equation, samaey2005gap, samaey2006patch}
provides a computational shortcut for approximating the homogenized solutions of systems when we only know the detailed problem.
The homogenized equation is \emph{assumed to exist, yet not available in closed form}.
The basic idea behind the equation-free approach is to exploit the \emph{implicitly assumed} existence of the macro-scale dynamics to speed up the micro-scale simulations and perform tasks like prediction, optimization or even bifurcation analysis at the large-scale level. In the homogenization case, the  micro-scale evolution is given by detailed, fine-scale PDE and the macro-scale description is the effective PDE. 
When a mathematically obtained homogenized PDE is available and accurate, we expect the equation-free numerics to accurately approximate the solutions of this homogenized PDE.

\begin{figure}
        \centerline{\includegraphics[scale=1.0]{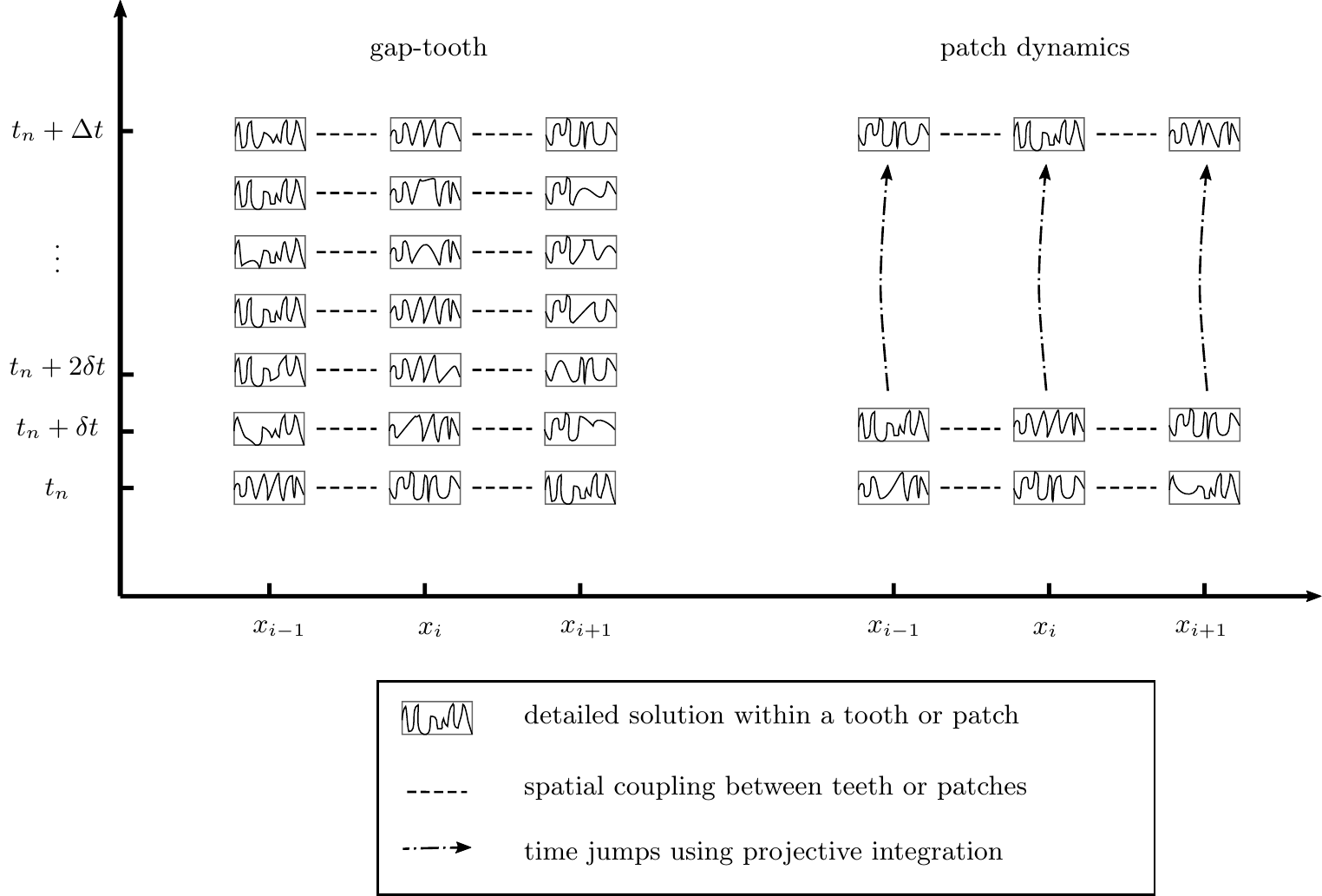}}
        \caption{ \textbf{Space-time grid of equation-free numerics.} In the gap-tooth scheme (left) the detailed problem is only simulated inside a sparse grid of teeth (boxes in the figure). The solution within adjacent teeth is coupled by implicitly utilizing the assumption that there exists a macro-scale variable which varies smoothly in the space. Similarly, in the patch dynamics, the simulation is carried out in sparse set of small boxes (called patches), but the smoothness of the macro-scale variable in time is leveraged to make the time grid sparse as well: the time derivative is estimated over the interval $\delta t$ and the solution is extrapolated over $\Delta t\gg\delta t$ using that estimate. }
        \label{fig_eqfree}
\end{figure}

Here, we utilize the equation-free approach for parsimonious generation of data: as explained below, the equation-free methodology guides us to simulate the detailed problem only in a fraction of space-time, to obtain training data for learning the effective (homogenized) equation, and hence, it substantially reduces the computational complexity of collecting training data.  In describing the equation-free approach for homogenization, we closely follow~\cite{samaey2009equation} and we refer the reader to~\cite{samaey2005gap, samaey2006patch, roberts2007general} for further analysis of this approach in the context of homogenization. A simplified representation of the equation-free approach is  visualized in \cref{fig_eqfree}.
 
A basic algorithm in the equation-free toolbox for spatially distributed systems is the \emph{gap-tooth} scheme. In this approach, we populate the spatial domain with sparsely distributed small boxes (subdomains, called \emph{teeth}) that only cover a small fraction of the spatial domain. The width of each tooth,~$h$, is chosen to be considerably larger than~$\epsilon$ but still small compared to the domain size. The distance between the centers of the teeth,~$\Delta x$, is chosen large compared to~$h$, so that the teeth occupy a small fraction of the full space.
We only simulate the detailed problem within the teeth, with appropriately prescribed ``coupling" boundary conditions (see below) and then define the gap-tooth solution~$u_{g}$ to be the detailed solution averaged over the tooth and reported at its midpoint. Let~$x_i$ be the position of the $i$th~tooth center; then
\begin{align}\label{eq_gpsol}
u^i_g := \frac{1}{h}\int_{x_i-h/2}^{x_i+h/2} u(\xi,t)d\xi,\qquad i=1,\ldots,N
\end{align}
is the gap-tooth prediction\slash order-parameter\slash amplitude for the $i$th~tooth, which approximates the discretized macro-scale variable, denoted by~$U$.
In many scenarios~\cite[e.g.][]{Bunder2020a, Roberts2019a, bunder2017good} we can prove that such a macro-scale field~$U$ exists and is governed by an explicitly unknown PDE (here, hopefully, the homogenized PDE), and that variations of~$U$ can be faithfully estimated using its value on the grid of such teeth centers (rather than a highly resolved grid designed to capture the detailed, $\epsilon$-scale material property variations).
Based on the assumption of the appropriateness of such a  macro-scale field behavior, artificial  boundary conditions for each tooth have been devised for the meaningfully coupled simulation of the micro-scale problem across all teeth~\cite[e.g.][]{Bunder2020a, roberts2007general}. These teeth-edge conditions constitute an important ingredient of the gap-tooth scheme, as they couple the dynamics of adjacent teeth and allow the global pattern of the solution to emerge from the simulations in (seemingly) separate teeth. One typically assumes that the macro-scale field~$U$ resembles a polynomial between teeth, that is, 
\begin{align}\label{eq_toothpoly}
U(x,t_n)\approx p_i^k(x,t_n), \qquad x\in[x_i-h/2,x_i+h/2]
\end{align}
with $p_i^k$ denoting a polynomial of even degree~$k$ within the $i$th~tooth. The coefficients of this polynomial are determined so that it results in the same box averages as the detailed solution within tooth~$i$ and $k/2$ teeth to its left and right (see~\cite{samaey2009equation} for the explicit formula of~$p_i^k$). Now, for simulating the detailed solution within the tooth~$i$ during the entire duration of the next time step, we may use, as the tooth-edge Neumann condition, the slope of the polynomial ~\cite{roberts2007general}, that is, 
\begin{align}\label{eq_toothbc}
\frac{1}{\epsilon}\int_{x_i+h/2-\epsilon}^{x_i+h/2+\epsilon}\p_\xi u(\xi,t)d\xi&=\p_x p_i^k\big|_{x_i+h/2},  \notag\\
\frac{1}{\epsilon}\int_{x_i-h/2-\epsilon}^{x_i-h/2+\epsilon}\p_\xi u(\xi,t)d\xi&=\p_x p_i^k\big|_{x_i-h/2} .
\end{align}
In this way, the detailed evolution in $i$th~tooth  is informed by the polynomial approximation of the macro-scale field near that tooth, which in turn, is informed by the value of the macro-scale field in the neighboring teeth. Further analysis of the gap-tooth scheme and a discussion of its convergence and accuracy can be found in~\cite{samaey2005gap, roberts2007general}. In particular, it was shown that the gap-tooth solution scheme leads to a finite-difference approximation of the homogenized PDE itself.

A more advanced algorithm in the equation-free toolkit is \emph{patch dynamics with buffers}\cite{samaey2009equation}. In this scheme, the simulations of the detailed problem is performed not only in a fraction of space, but also in a fraction of time. This is made possible through the use of (coarse) projective integration \cite{gear2003projective,gear2002coarse}:  the gap-tooth simulations of the detailed problem are carried out over a few short time steps, and the result is used to compute an estimate of the time-derivative for the macro-scale variable,~$\p_t U$. This estimate is then used to take a large time step in the simulation of~$U$. Taking large time steps for~$U$ is plausible, because we have assumed the existence of a macro-scale closure for~$U$, which justifies using coarser grids in time compared to those necessary for the accurate simulation of the detailed problem.
However, and similar to the gap-tooth scheme, the macro-scale closure is not available, and thus cannot be explicitly used in the algorithm. 
Instead, the quantities required for scientific computation (like the time derivatives here) are \emph{estimated} from judiciously designed partial simulations of the detailed problem.

The grid setup for patch dynamics with buffers is similar to the gap-tooth scheme, except that we extend each patch to width $H> h$, with a buffer zone around each tooth. 
The buffers provide an alternative way of imposing meaningful coupling boundary conditions across neighboring teeth~\cite{bunder2017good}. 
Detailed simulations \emph{are also performed} in the buffer zones, yet the patch dynamics solution reports averages over the tooth itself. 
The role of the buffer zone is thus to ``protect" the core solution within the tooth from the disturbances caused by the boundary conditions on the buffers. 
Now the teeth exchange information through the process of \emph{lifting} defined below.

To round off our brief outline of the patch dynamics scheme, we define two more notions in multi-scale computation.
Firstly, the \emph{restriction operator} is a function that maps the detailed solution within patches to the macro-scale variable. 
Here, the restriction operator is simply averaging over each tooth in~\eqref{eq_gpsol}. 
Secondly, the \emph{lifting operator} is a (one-to-many) mapping from the macro-scale field to a consistent detailed solution profile. The lifting operator is not unique, since there are many detailed profiles that have the same average. A standard choice for lifting in homogenization is a polynomial expansion similar to \cref{eq_toothpoly} within the $i$th~tooth \emph{and its surrounding buffer zone}. In the explicit form, this polynomial can be written as
\begin{align}\label{eq_lift}
u_i(x)=\sum_{k=0}^dD_i^k(U)\frac{(x-x_i)^k}{k!},\quad x\in\left[x_i-\tfrac{H}{2},x_i+\tfrac{H}{2}\right],
\end{align}
where~$D_i^k(U)$ is the finite-difference approximation for the $k$th~spatial derivative at the $i$th~tooth.

A time step of the patch dynamics simulation is performed as follows: given the current state of the macro-state field $U_i^t,~i=1,\ldots,N$, the state of the detailed solution within  each tooth and its buffer is constructed using the lifting in~\eqref{eq_lift}. Then the detailed equation is marched forward in the teeth and their buffers. The new macro-state $U_i^{t+\delta t},~i=1,\ldots,N$, is computed by averaging in~\eqref{eq_gpsol}. 
Then, the time-derivative of the macro-state in each tooth is estimated as
\begin{align}\label{eq_dUdt}
\frac{d}{dt}U_i(t)\approx\frac{U_i^{t+\delta t}-U_i^t}{\delta t}.
\end{align}
Next we perform a projective integration step: we use the above estimate within any time-stepping method (like forward Euler or Runge--Kutta), to march~$U$ forward in (long) time, and in particular, use a large time step $\Delta t\gg\delta t$ to approximate the new macro-state $U_i^{t+\Delta t},~i=1,2,\ldots,N$.   The patch dynamics scheme, effectively, requires us to simulate the detailed problem only in a fraction of space (i.e., teeth and buffer) and a fraction of time (i.e., brief sequences of short time steps for computing the gap-tooth estimate~\eqref{eq_dUdt}). 
Samaey et al. \cite{samaey2006patch} discussed in detail this patch dynamics algorithm for homogenization including its convergence and choice of boundary conditions.

\section{Learning of effective PDEs from data}\label{sec_learning}
The key contribution of this work is exploiting the equation-free approach for parsimonious generation of training data in learning effective homogenized PDEs. In general, although homogenized equations may exist for heterogeneous materials, there may be no systematic way to analytically discover them (for example in the case of non-periodic heterogeneity~\cite{pavliotis2008multiscale}). In such cases, one seems restricted to fully simulating the detailed problem,  and then \emph{performing coarse-grained observations of it} (e.g., via averaging or low-pass filtering) to obtain training data.

\begin{figure}
        \centerline{\includegraphics[scale=1.0]{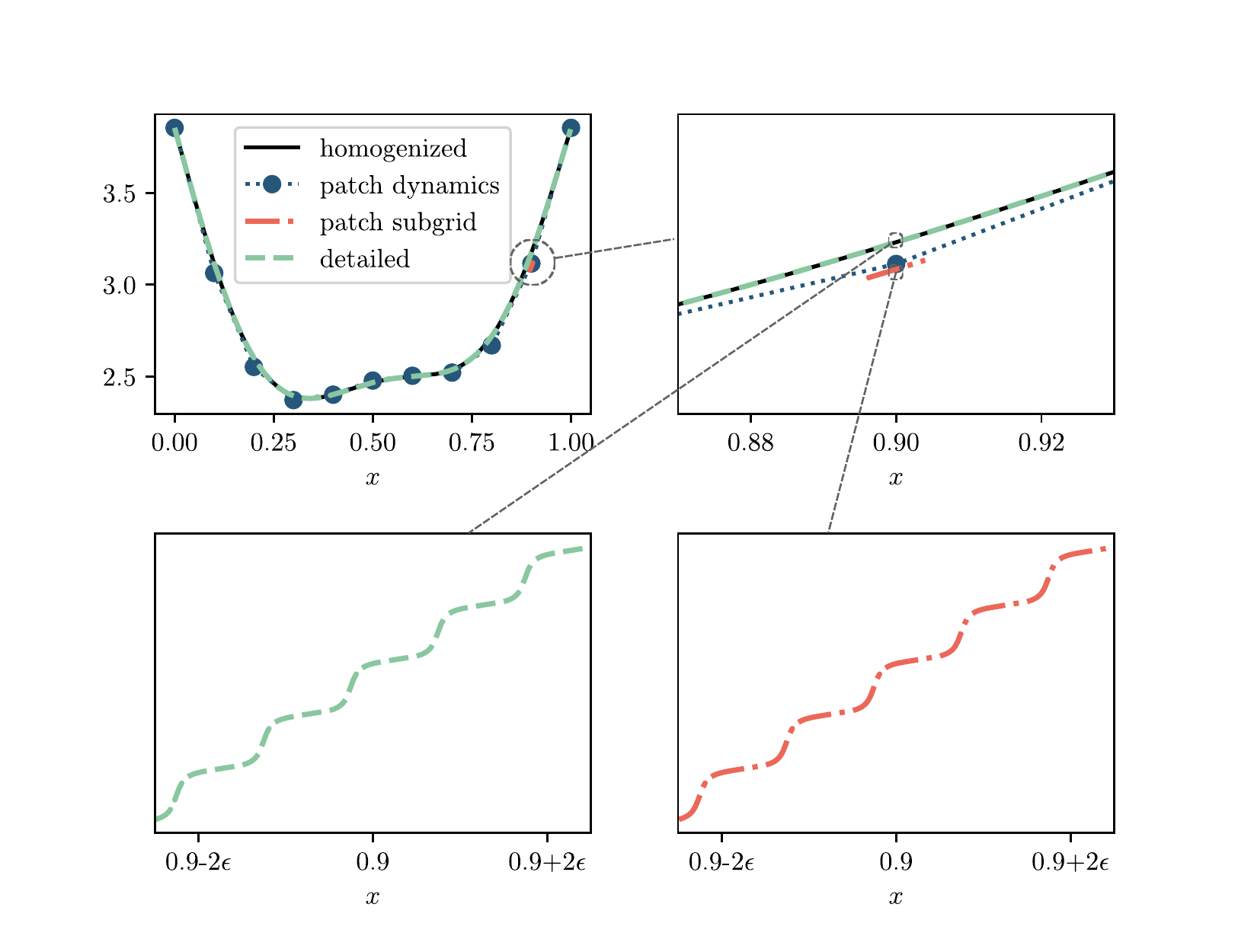}}
        \caption{ \textbf{A sample snapshot from different computational implementations of the one-dimensional heterogeneous diffusion problem  in \cref{sec_result1}.} The top left panel demonstrates the agreement between the fully resolved, the homogenized, and the patch dynamics solution observed at large-scales. The bottom two panels show magnifications revealing the $\epsilon$-scale variations of the solution for the fully resolved and the patch dynamics simulation ($\epsilon=10^{-5}$). Instead of the fully resolved simulation, we use a patch dynamics sparse simulation ($ 10^{4}\times$ smaller space-time grid) to collect data for learning the effective dynamics. This particular snapshot is produced from an initial condition of the form~\eqref{eq_IC1d} after running the simulations for 0.05 time units.}
        \label{fig_detailed}
\end{figure}

Here, we use the equation-free approach as the alternative, as it only requires simulation of the detailed problem only in a fraction of the total space and time domain. \Cref{fig_detailed} shows sample snapshots of data produced in an example by  solving the homogenized PDE, the detailed PDE, and the  patch dynamics approximation of the effective PDE for the one-dimensional diffusion problem in  \cref{sec_result1}. 
In this example, the patch dynamics allows us to  simulate the detailed problem in less than 8\% of the space and~$0.1\%$ of time domain, effectively reducing the space-time grid size by a factor of~$10^{4}$. 
We now use the data so produced to train neural nets that learn the time-evolution (i.e. the effective operator) of the effective (homogenized) system. To validate the approach, we evaluate the performance of the trained neural nets via comparison with the homogenized equations which \emph{are known in closed form} for our test examples. We remind the reader that such closures are, in general, not available.

We consider two neural network learning approaches for discovering the homogenized PDE from the data. In the first architecture, we try to discover the \emph{functional form} of the law of the PDE, i.e., the relation between spatial and temporal derivatives of the coarse field  at every macro-scale space-time point:
\begin{equation}\label{eq_PDE1}
\p_t U=F\big(U,\p_x U,\p_{xx} U, \ldots).
\end{equation}
The particular combination of derivatives required for successful learning  of the homogenized equation may be unknown \emph{a priori}. In previous work~\cite{lee2020coarse}, we successfully utilized learning approaches such as diffusion maps and Gaussian process regression to identify a minimal subset of spatial derivatives that suffices to learn an effective model.
Assuming that we have included a sufficient set of spatial derivative terms as input,  the learning problem amounts to regressing the function~$F$ on its arguments (i.e., learning the right-hand-side differential operator~\cite{chen1995universal, lu2019deeponet}) from data.

In the second approach, we model the law of the PDE directly in a \emph{discretized} form: as a relation between the time derivative of the field at a point of interest and the values of the field at a number of points in the neighborhood of the point of interest. This can be thought of as a form of spatial discretization of the effective PDE for the macro-scale variable 
 \begin{equation}\label{eq_PDE2}
\frac{d}{dt} U(x_j)=G\big(U(x_j),U(x_{j-1}),U(x_{j+1}), \ldots).
\end{equation}
where $x_j$'s denote points on a spatial grid around the point of interest ordered by the index~$j$. The stencil for this learning problem should be large enough to contain sufficient information for successfully approximating the local~$dU/dt$. 

The above two architectures appear very similar: they employ the local profile in the neighborhood of the grid point of interest either (a)~in the form of observed\slash estimated spatial derivative values \emph{at that point}, or (b)~in the form of a grid of neighboring field values, in order to approximate time evolution.
Both approaches are capable of learning PDEs with spatially homogeneous character with various terms such as diffusion, reaction, linear or nonlinear advection and source or sink terms.
The two approaches can be formally related to each other when the observed dynamics (possibly after fast initial transients) lie on a low-dimensional manifold, for example on a finite-dimensional inertial manifold~\cite{jolly1990approximate, constantin2012integral}. Theoretical justification for the second approach can also be provided through the center-manifold-based work in~\cite{roberts2001holistic} towards the ``holistic discretization", which  has been taken advantage of in enhancing patch dynamics simulations~\cite{bunder2016accuracy, samaey2009equation}.

A practical advantage of the second approach is that it avoids the ``off-line" estimation  of spatial derivatives, which could be sensitive to noise and computational directional preferences (e.g., upwinding for advection problems); the obvious price to pay---along with the loss of translational invariance due to discretization---is the dependence of the learned network parameters on the particularly selected grid. Therefore, the learned model has to be retrained or adjusted for new spatial grids in the second architecture. This second approach has already been utilized in learning data-driven discretization schemes~\cite[e.g.][]{bar2019learning}.

\begin{figure}
\centering
\subfloat[][]{\includegraphics[scale=1]{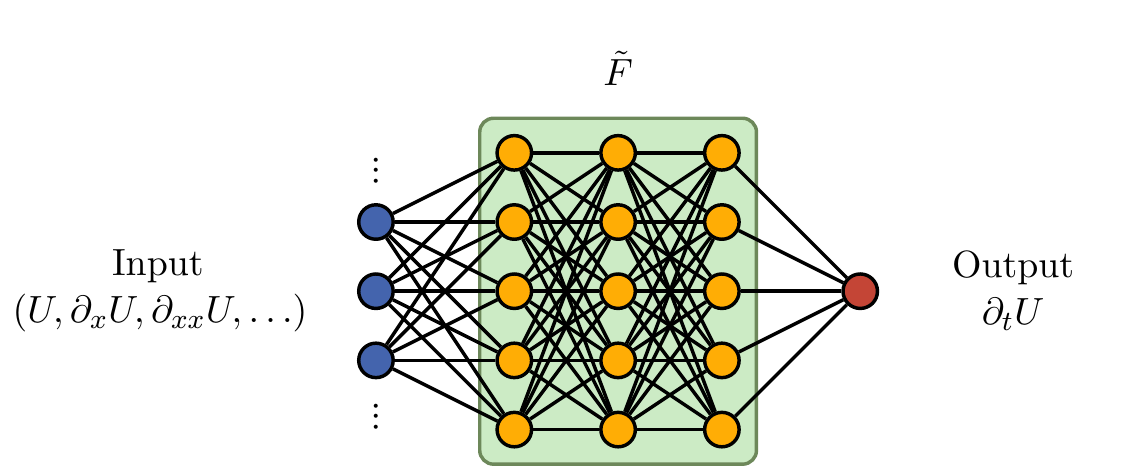}} 
 
\subfloat[][]{\includegraphics[scale=1]{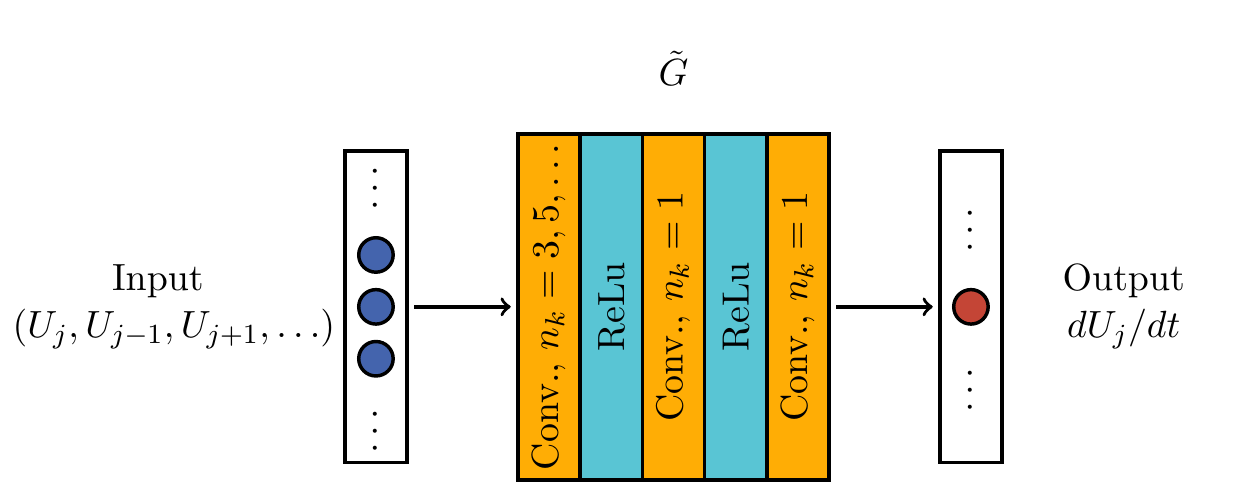}}

\caption{\textbf{Neural network architectures for approximate PDE discovery.} (a)~A neural net architecture for learning the functional form of the  PDE differential operator: there are three fully connected layers with 32 nodes and the first two layers have ReLU nonlinearity. (b)~A neural net architecture for learning the PDE in its discretized form: there are three convolutional layers with 32 filters and ReLU activations in between. The last two convolution layers have kernels of size one and, in effect, act as fully connected layers.}
\label{fig_nnarch}
\end{figure}

We use the neural network architectures shown in \cref{fig_nnarch}(a) and (b), respectively, to regress~$F$ and~$G$. The first network consists of three fully connected layers with a Rectifying Linear Unit~(ReLU) as the activation function in the first two layers. 
Recall that the action of a fully connected linear layer can be written as
\begin{align}
z^{out}=Az^{in} + b, \quad b,z^{out}\in \R^n,~z^{in}\in\R^m,~A\in{\R^{n\times m}}.
\end{align}
On the other hand the action of ReLU layer is given by
\begin{align}
z^{out}=\max(0,z^{in}).
\end{align}
Therefore our network represents a function of the form
\begin{align}
z^{out}=A^{(3)}\max\big(A^{(2)}\max\left(A^{(1)}z^{in} + b^{(1)},0\right) + b^{(2)},0\big)+b^{(3)}
\end{align}
where the superscripts are the index of the fully connected layer. Entries of the matrices $A^{(i)},~i=1,2,3$ and the bias vectors  $b^{(i)}, ~i=1,2,3$ are the parameters of the network which are to be optimized in the learning process.

The second network consists of three convolutional layers, with ReLU activation layers in between.Let us recall the mathematical description of a (linear) convolution layer action. Assume each layer consists of $n_k$ kernels of size $2n+1$. Let $C_{lpq}$ denote the $q$-th kernel in the layer. Then the action of the layer on an input with $p$ channels is given by
 \begin{align}
z^{out}_{iq}=\sum_p\sum_{j=-n}^{n}C_{jpq} z^{in}_{i+jp} + b_{q},\quad q=1,\ldots,n_k,
\end{align}
which, in the shorthand form is written as 
 \begin{align}
z^{out}=C*z^{in}.
\end{align}
As such, our network represents the function
 \begin{align}
z^{out}=C^{(3)}*\big(C^{(2)}*\max(C^{(1)}*z^{in} + b^{(1)},0)+b^{(2)},0\big)+b^{(3)}
\end{align}
where the superscripts are the index of the convolutional layer in the network. The entries of the tensors $C^{(i)},~i=1,2,3$ and vectors $b^{(i)},~i=1,2,3$ are to be learned from data.
The first convolutional layer, with a sufficiently large kernel size, effectively acts as a finite-difference model of the spatial derivatives. Here, we set the kernel size of the next two layers to be one, which makes them, in effect,  fully connected layers.
This general type of shallow architecture has been demonstrated to suffice for capturing the nonlinearity of PDEs arising in fluid mechanics~\cite{bar2019learning}. 

\clearpage
\section{Results}\label{sec_results}
\subsection{One-dimensional heterogeneous diffusion}\label{sec_result1}

In the first example, we consider the one-dimensional heterogeneous diffusion problem from~\cite{samaey2006patch}. The detailed evolution is governed by the PDE
\begin{align}\label{eq_detailed1d}
\p_t u(x,t)=\p_x\left[ a\big(\frac{x}{\epsilon}\big) \p_x u(x,t)\right],\quad a\big(\frac{x}{\epsilon}\big)=1.1+\sin\big(2\pi \frac{x}{\epsilon}\big)
\end{align}
on the domain~$[0,1]$ with $\epsilon =10^{-5}$ and time-independent Dirichlet boundary conditions of
\begin{align}
u(0,t)=u(0,0),\quad u(1,t)=u(1,0).
\end{align}
The corresponding homogenized equation is known to be
\begin{align}\label{eq_hom1d}
\p_t u_h(x,t)=\p_x\left[ a^*\p_x u_h(x,t)\right],\quad a^*\approx 0.45825686,
\end{align}
with similar boundary conditions. 

To generate the training data, we use the patch dynamics scheme for \cref{eq_detailed1d}  using 10 uniformly spaced teeth in the spatial domain ($\Delta x=0.1$).
Each tooth has a width of $h=10^{-4}$ and it is surrounded by a buffer zone of width $H=8 \times 10^{-3}$. The detailed simulation within each tooth and its buffer is carried out using a central finite-difference scheme with spatial resolution of $\delta x=10^{-7}$ and time step size of $\delta t=10^{-6}$ with Dirichlet boundary conditions imposed on the buffer zone boundaries. 
The macroscopic time step for advancing the macro-state variable is $\Delta t=10^{-3}$. Therefore, we are only solving the detailed problem in 8\% of the spatial domain and~$0.1\%$ of the time domain. \cref{fig_detailed} shows sample snapshots of~$u$ for this problem arising during the computation using various alternative methods.

We simulate 8 trajectories of the system with initial conditions of the form 
\begin{align}\label{eq_IC1d}
u(x,0)=\sum_{j=1}^{20}a_j\sin(2\pi l_j x + \phi_j)
\end{align}
where $a_j,~l_j$ and~$\phi_j$ are randomly drawn from uniform distributions on~$[-1,1]$, $[0,4]$ and~$[0,2\pi]$, respectively. We simulate each trajectory over the time interval of~$[0,1]$ and record snapshots of~$u$ and~$\p_t u$ every $\delta t_s=10^{-3}$.
We train both neural net architectures using the ADAM optimizer~\cite{kingma2014adam} with batch sizes of 64 and learning rate of~$10^{-3}$. For the optimization objective, we use the Mean Squared Error~(MSE) of the network output (i.e.,~$\p_t u$). In this example, we use a second-order finite-difference method with stencil size of three  to estimate the spatial derivatives required as inputs by the first architecture, and a corresponding kernel size of three in the first convolutional layer of the second architecture.

\begin{figure}
        \centerline{\includegraphics[width=1 \textwidth]{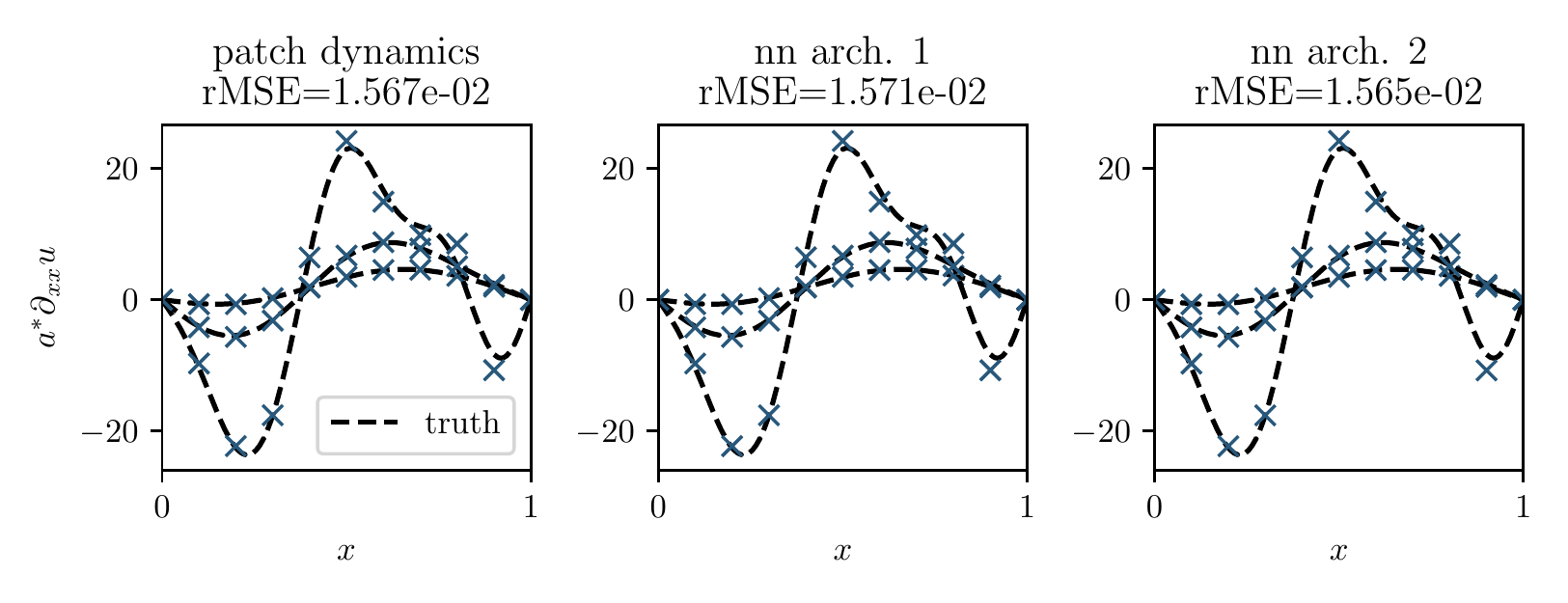}}
        \caption{ \textbf{Prediction of the RHS of an effective PDE via patch dynamics and neural nets}. ``Truth" here denotes the solution of the homogenized equation on a fine grid.}
        \label{fig_dudt1d}
\end{figure}

To validate the learned models, we compare their performance to that of the known homogenized model in two test trajectories that have initial conditions of the form~\eqref{eq_IC1d} but have not been used in the training. The (converged) homogenized solution to which we compare is computed using a second-order finite-difference discretization of \cref{eq_hom1d} on a grid size of $dx=5\times 10^{-3}$ and \verb1lsode1 time stepper. As the error measure, we use the relative Mean Squared Error (rMSE) between the two snapshots, defined as
\begin{align}\label{eq_rMSE}
\text{rMSE}(v)=\frac{\sum_j\|v^j - v^j_h\|^2}{\sum_j\| v^j_h- \overline{v_h}\|^2 }
\end{align}
where $v^j$ denotes the $j$th~snapshot predicted by the neural net  and $v^j_h$ is the associated snapshot from the homogenized model. As shown in \cref{fig_dudt1d}, both architectures are quite successful in accurately representing the RHS of the homogenized PDE, even though the training data is generated on a much coarser grid than the fully resolved data.

\begin{figure}[!h]
        \centerline{\includegraphics[width=1 \textwidth]{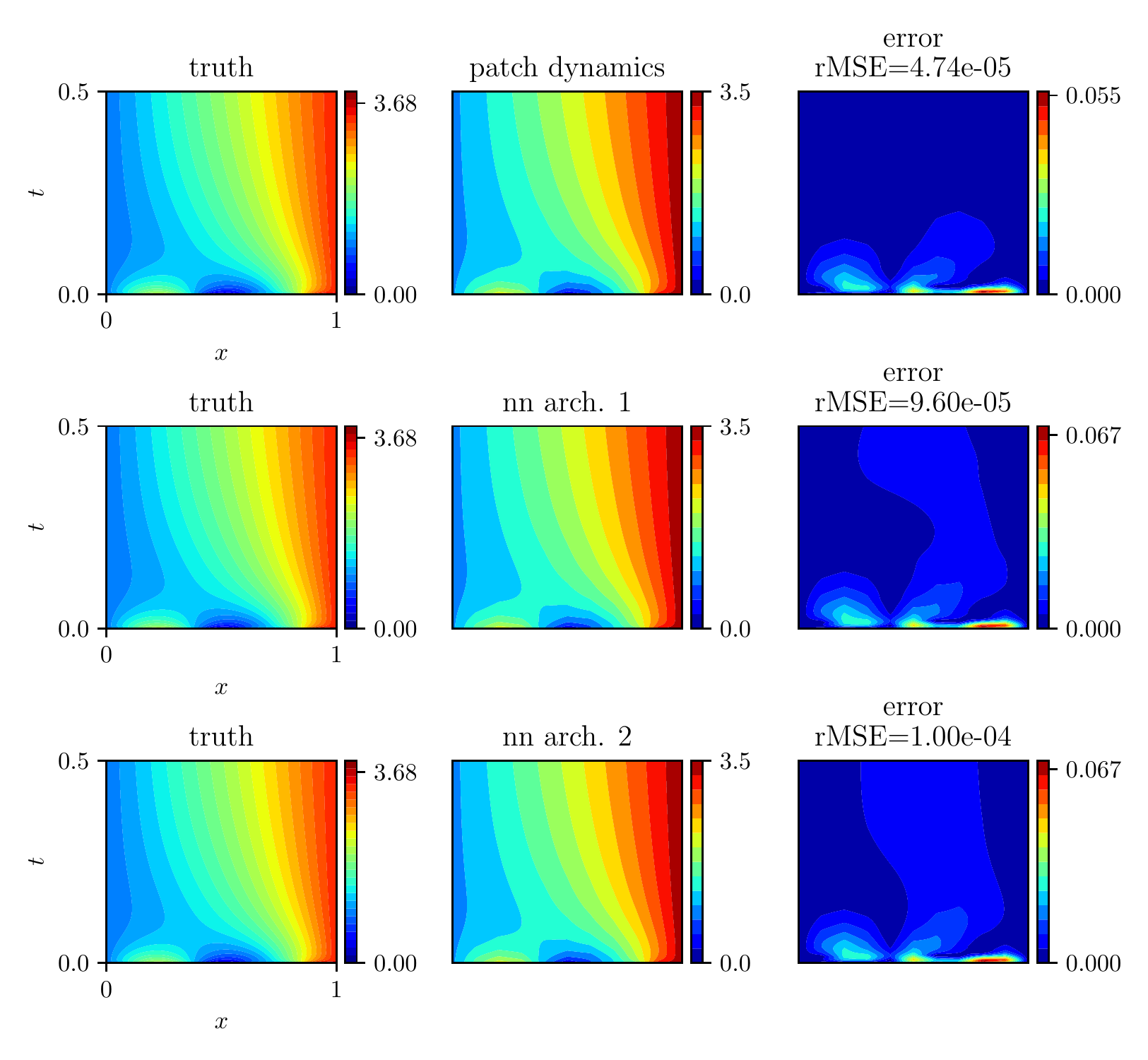}}
        \caption{\footnotesize \textbf{Integration of a trajectory from an unseen initial condition via patch dynamics and neural nets}. ``Truth" here denotes the solution of the homogenized equation on a sufficiently resolved grid.}
        \label{fig_traj1d}
\end{figure}

Next, we use the learned models to simulate a trajectory with a new initial condition from~\eqref{eq_IC1d}. We use the \verb1BDF1 method~\cite{shampine1997matlab} as the time stepper on the same macroscopic grid. Comparison to the solution of the homogenized PDE (right column of \cref{fig_traj1d}), suggests that the learned models  exhibit, in addition to accuracy,  the stability properties required for time integration of the effective equations.

\clearpage
\subsection{Diffusion in two dimensions with lattice heterogeneity}

In this example~\cite{bunder2017good,Bunder2020a},  we consider diffusion on a two-dimensional periodic domain  with lattice-type micro-scale heterogeneity. We define a two-dimensional lattice, indexed by~$(i,j)$, with the same spacing~$h$ in both directions on the bi-periodic domain~$[0,2\pi)^2$. The evolution of the micro-scale field~$u_{i,j}$ is given by
\begin{align}\label{eq_2dhet}
\dot{u}_{i,j}=\big[\K^{x}_{i,j} (u_{i+1,j}-u_{i,j})+\K^{x}_{i-1,j}(u_{i-1,j}-u_{i,j})
+\K^{y}_{i,j}(u_{i,j+1}-u_{i,j}) + \K^y_{i,j-1}(u_{i,j-1}-u_{i,j})\big]/h^2
\end{align}
where $\K^{x}_{i-1,j}$ is the diffusivity at  half lattice point~$(i+1/2,j)$ and $\K^{y}_{i-1,j}$ is the diffusivity at half lattice point~$(i,j+1/2)$. We assume that the heterogeneous diffusivity is 3-periodic in both directions and, in particular, given as
\begin{align}\label{eq_chetr}
\K^x=
\begin{bmatrix}
1.0566 & 0.6668 & 1.1568\\ 6.5894 & 0.8683  &  2.4174 \\ 0.9473  & 1.1407 &  1.6610
\end{bmatrix},
\quad
\K^y=
\begin{bmatrix}
3.6355  & 0.4470 &   2.3896\\
    0.8628  &  4.8558  &  0.2833\\
    4.5025  &  1.5865  &  0.5679
\end{bmatrix}.
\end{align}
It is known \cite{bunder2017good} that there exists a homogenized equation, namely,
\begin{align}\label{eq_2dhom}
\p_t u_h = 1.2644 \p_{xx} u_h + 1.3398 \p_{yy} u_h
\end{align}
which provides an approximation to the evolution of~$u$ at large length scales.

We use the gap-tooth scheme for parsimoniously simulating the micro-scale system. We use a uniformly spaced $16\times 16$ grid of square teeth, each one with sides of $H=4h=0.0524$. These spatial patches approximately cover~1\% of the space. 
We solve \cref{eq_2dhet} for the micro-scaled grid points that are contained within our spatial patches under the appropriate boundary conditions for each patch. The macro-scale solution is then computed by averaging~$u$ on the interior points of the patch. Here, the simulation is performed for all times (no projective integration ``jumps" are used); in equation-free terms this is a gap-tooth scheme, and the savings come only from the spatial sparsity of the simulation. 
The spatial patch boundary conditions are carefully chosen to couple the macro-scale solution on neighboring patches and approximately retrieve the global coarse pattern that would appear from simulating the fully resolved micro-scale system, as detailed elsewhere~\cite{bunder2017good}.

We simulate 100 trajectories where the initial conditions are given by 
\begin{align}\label{eq_IC2d}
u(x,0)=\sum_{j=1}^{10}a_j\sin( l^x_j x + \phi^x_j)\sin( l^y_j y + \phi^y_j)
\end{align}
and the parameters $a_j\in [-1,1]$, $(l^x_j,l^y_j)\in \{1,2,3,4,5\}^2$, and $(\phi^x_j,\phi_j^y)\in[0,2\pi]^2$ are drawn randomly from uniform distributions on their domain. We simulate the trajectories over $t\in[0,1]$ and collect 100 snapshots of each at sampling intervals of~$0.01$. 
We use data from 85 trajectories (i.e., 8500 snapshots) to train our networks. The training parameters are chosen similar to the one-dimensional problem. In this example, we use FFT to estimate the spatial derivatives for the first neural net architecture, thus improving the estimates of the input to the neural net.

\begin{figure}
        \centerline{\includegraphics[width=1 \textwidth]{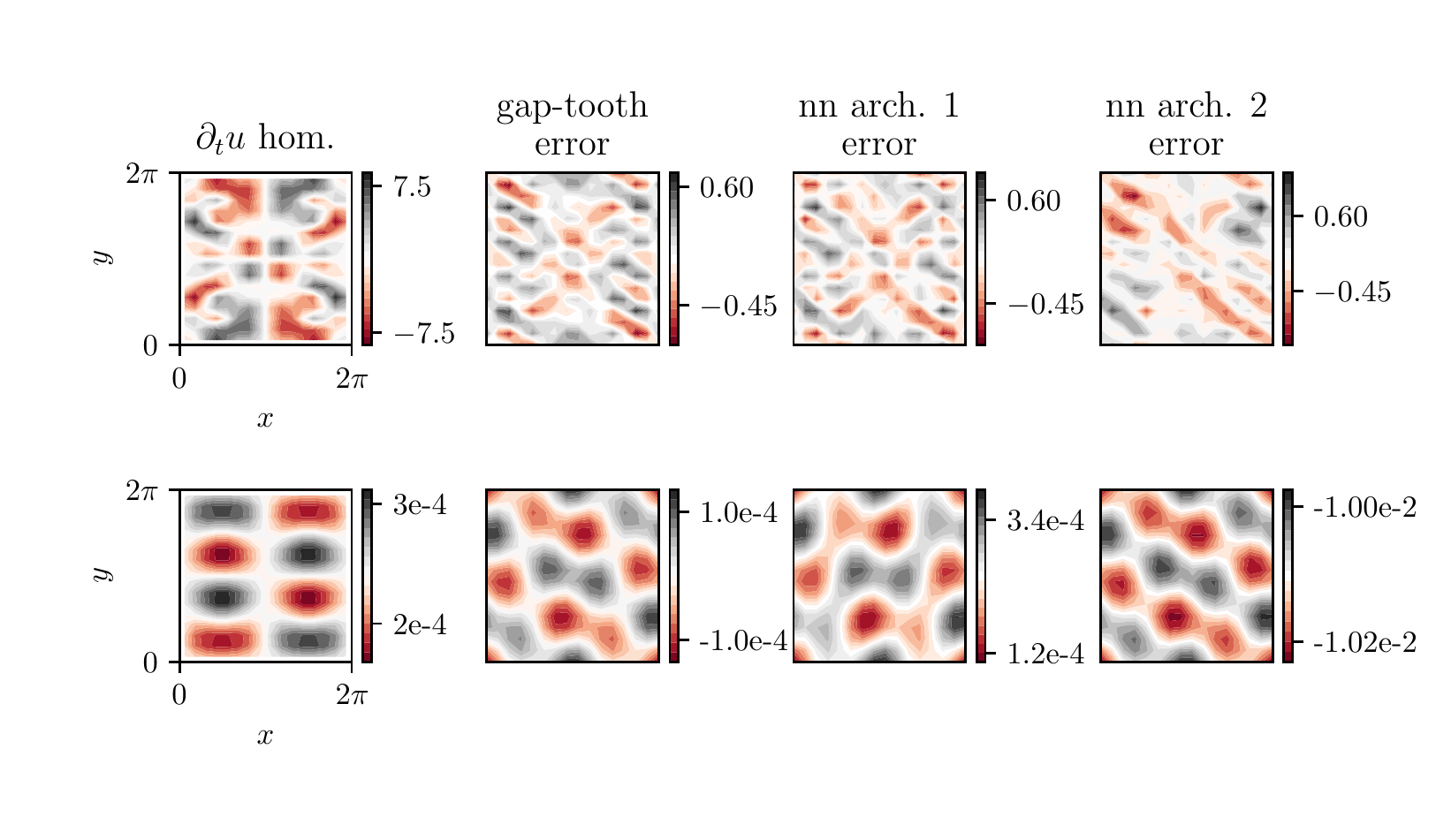}}
        \caption{ \textbf{Performance of models in estimating 2D PDE's right-hand-side.} The relative error of prediction is small for snapshots with large magnitude of~$\p_t u$ (e.g., first row) and large for  small magnitudes of~$\p_t u$ (e.g., bottom row). }
        \label{fig_2dut}
\end{figure}

In this example our scheme, trained with gap-tooth simulation data averaged over a sparse grid of teeth, is less successful in providing an accurate closed evolution equation. 
Our neural networks (with various depths) cannot drive the MSE significantly below a few percent.
We believe that this shortcoming can be rationalized: 
\Cref{fig_2dut} shows sample snapshots of~$\partial_t u$ from the homogenized solution, which we take as the truth here, and the prediction error of our best trained neural networks.  All the trajectories in this example converge to steady state where $\p_t u=0$. We observe that the relative prediction error becomes typically larger for snapshots with smaller magnitude of~$\p_t u$. A rationalization of this observation is that the choice of loss function (e.g., MSE in this work) emphasizes accurate prediction \emph{for large-magnitude snapshots}, and that for the small-amplitude data used here, it might be expected that the relative error would then be large. Therefore, the learned estimation of the time-derivative near the steady state is less accurate than estimates for states away from steady state.


\begin{figure}
        \centerline{\includegraphics[width=1\textwidth]{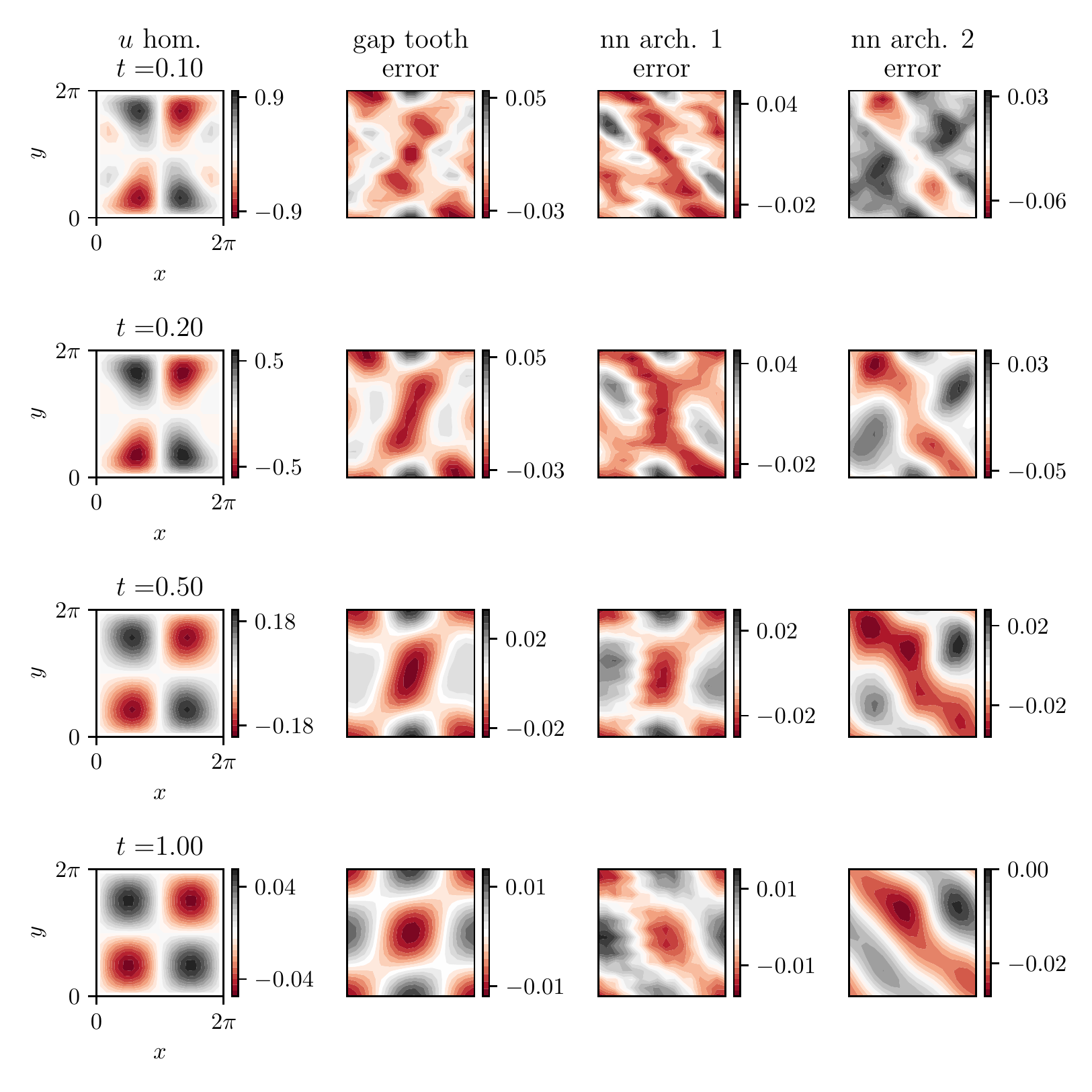}}
        \caption{ \textbf{Performance of learned models in integrating a test trajectory.} The neural net models (trained via gap-tooth data) follow the homogenized solution closely at the beginning but yet lag behind the final decay to steady state. }
        \label{fig_2du}
\end{figure}

\begin{figure}
\centering
\subfloat[][]{\includegraphics[scale=1]{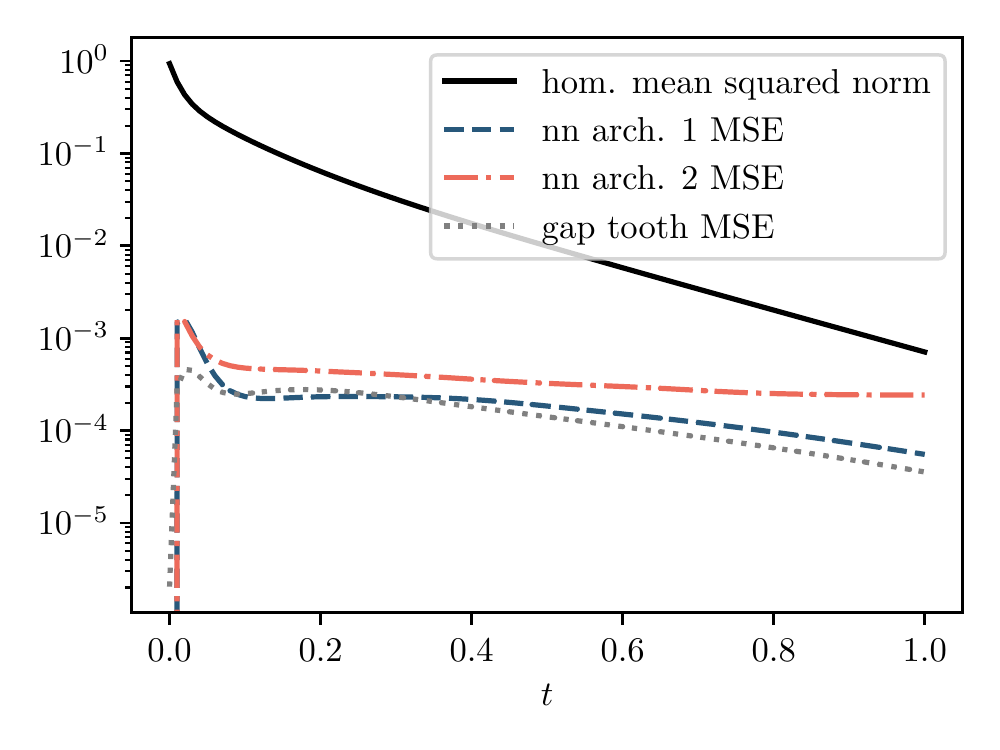}} 

\subfloat[][]{\includegraphics[scale=1]{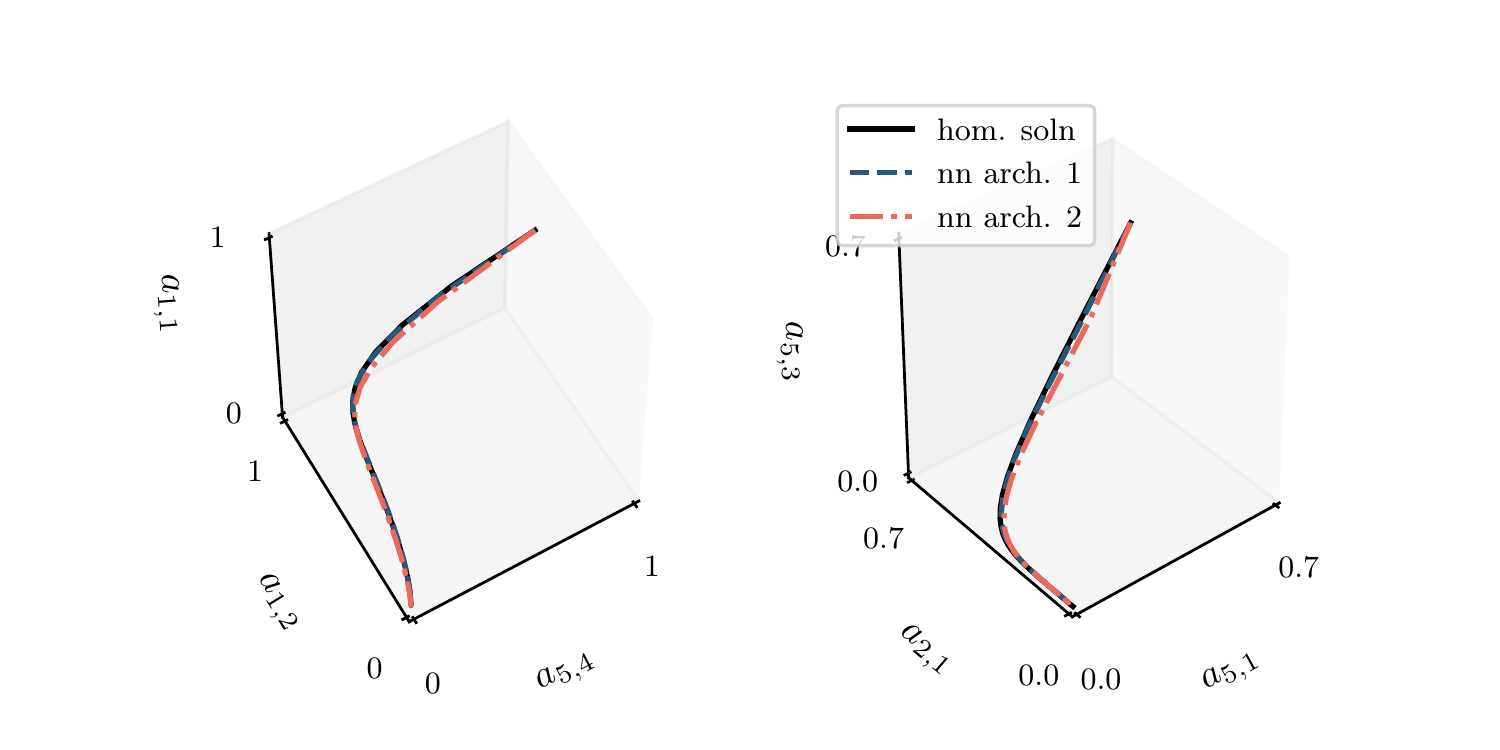}}
\caption{\textbf{Qualitative behavior of learned models} (a)~The mean squared error~(MSE) of test trajectories and (b)~a state space viewpoint showing the six largest Fourier amplitudes of the solution trajectory.}
\label{fig_ss}
\end{figure}

\Cref{fig_2du} shows sample snapshots of a trajectory integrated using the PDEs discovered by the neural nets and their comparison with the corresponding homogenized solution. The neural net solutions closely follow the homogenized solution early in the integration interval (at large solution amplitudes) and their error is comparable to that of the training data (i.e., gap-tooth versus homogenized). However, as the true trajectory approaches the steady state, and the solution amplitude decays, the relative error of the learned models grows---the time derivative is less accurately estimated in that region, and that error accumulates in the trajectory over time (\cref{fig_ss}(a)). Interestingly, however, \emph{the qualitative state-space behavior of the discovered models is visibly accurately recovered, and close to truth:} \cref{fig_ss}(b) shows a projection of the solution trajectories from the homogenized solution and from the neural-nets-identified PDE in the solution state space.  This projection is realized by observing the six largest coefficients of the Fourier expansion for the solutions. The trajectories from the learned models, regardless of time parametrization, are visibly close to the homogenized solution trajectory, and only start to deviate in a small neighborhood of the origin. We then credibly argue that the overall identification process is successful: the trajectories of the learned models quantitatively shadow the true ones to the neighborhood of the final steady state, and even there they have the right long term\slash stability characteristics.

\clearpage
\section{Conclusion}\label{sec_conclusion}
Machine learning techniques have proven successful in learning PDEs that approximate the observed evolution of complex systems from spatiotemporal data. Despite this initial success, an outstanding challenge for identifying models of the effective dynamics of multiscale\slash multiphysics systems is the computational burden of generating the training data from detailed computational models embodying our best understanding of the fine-scale physics. 
In this work, we proposed linking established multiscale scientific computing algorithms (and, in particular, using the equation-free approach, with modern neural network algorithms for learning effective PDEs. We illustrated the combination of these two ``computational technologies"  for learning homogenized equations for materials with micro-scale property heterogeneities, enabled by parsimonious training data collection. The advantage of our approach lies precisely in the efficient compilation of training data: the equation-free approach (and other similar approaches~\cite[e.g.][]{weinan2007heterogeneous}) requires simulation of the detailed problem in only a small fraction of the space-time domain and dramatically reduces the amount of computational effort for collecting the requisite training data. We demonstrated our approach in the case of heterogeneous diffusion and validated our results using the analytically known homogenized solutions. Our approach can be readily extended to problems with more complex material behavior (non-periodic heterogeneity and drift terms) where the there are no established methods to derive closed-form solutions (assuming such equations exist).

Throughout this paper, we assumed the appropriate macro-scale field variable is known. Indeed, guided by the classical homogenization theory, we took the target macro-scale variable (i.e., homogenized equation field) to be a low-pass filtered version of the micro-scale variable (the fully resolved solution), and used the averaging-over-tooth operation within the equation-free approach to transform fine-scale data to training data for the macroscopic model (what in equation-free terminology is called ``the restriction operator").  An important research direction involves the data-driven identification of the right observables, i.e. the variables in terms of which the PDE to be learned should be formulated. Whether through manifold learning techniques, like diffusion maps (e.g., in the ``variable free" computational efforts)~\cite{erban2007variable, chiavazzo2017intrinsic, liu2015equation}
 or through deep-learning tools like (possibly variational) auto-encoders, and Wasserstein GANs~\cite{arjovsky2017wasserstein, rico1992discrete, kingma2013auto},
 this possibility will, we expect, greatly enhance the range of phenomena for which useful data-driven models can be derived.
 Physically informing such models~\cite[e.g.][]{raissi2019physics} and providing explainable descriptions of their predictions is a complementary, but equally important, challenge.
One can even contemplate the case of ``emergent PDEs", where not only the dependent variables and the operators, but even the independent variables themselves (the right ``space" and even the 
right ``time" for the PDE) are obtained in a data-driven manner~\cite{kemeth2018emergent, thiem2020emergent}.

\section*{Acknowledgements}
The authors declare that they have no conflict of interest in publication of this paper. This work was partially supported by a US ARO MURI (through UCSB), by  DARPA, and by the Australian Research Council grant DP200103097. Discussions with Dr.~T. Bertalan are gratefully acknowledged.

\section*{Data and source code }
A python Implementation of our framework and the data for producing the figures in this paper is available at \url{https://github.com/arbabiha/homogenization_via_ML}. 

\printbibliography
\end{document}